\documentclass[graybox]{svmult}

\usepackage{mathptmx}       
\usepackage{helvet}         
\usepackage{courier}        
\usepackage{type1cm}        
\usepackage{makeidx}         
\usepackage{graphicx}        
\usepackage{multicol}        
\usepackage[bottom]{footmisc}

\makeindex             

\begin{document}

\title*{Constructive martingale representation in functional It\^o calculus: a local martingale extension}
\titlerunning{Constructive representation of local martingales} 
\author{Kristoffer Lindensj\"o}
\institute{Kristoffer Lindensj\"o \at Department of Mathematics, Stockholm University, SE-106 91 Stockholm, Sweden 
\email{kristoffer.lindensjo@math.su.se}
}
%
%
\maketitle

\abstract{The constructive martingale representation theorem of functional It\^o calculus is extended, from the space of square integrable martingales, to the space of local martingales. The setting is that of an augmented filtration generated by a Wiener process.}

\vspace{8mm}

\noindent {\footnotesize First published in:} \textcolor{blue}{\footnotesize chapter 9 in S. Silvestrov et al. (eds.), 'Stochastic Processes and Applications', Springer Proceedings in Mathematics \& Statistics 271, 165--172, 2018.}
{\footnotesize DOI:}\textcolor{blue}{\footnotesize 10.1007/978-3-030-02825-1}. {\footnotesize Unauthorized reproduction prohibited.}

\section{Introduction} \label{introLINDENSJO}
Consider a complete probability space $(\Omega,{\mathcal F},\bbbp)$ on which lives an $n$-dimensional Wiener process $W$. Let $\underline{\mathcal{F}} =  (\mathcal{F}_t)_{0\leq t \leq T}$ denote the augmentation under $\bbbp$ of the filtration generated by $W$ until the constant terminal time $T<\infty$. 
One of the main results of It\^o calculus is the martingale representation theorem which in the present setting is as follows: \emph{Let $M$ be a RCLL local martingale relative to $(\bbbp,\underline{\mathcal F})$, then there exists a progressively measurable $n$-dimensional process $\varphi$ such that}

\begin{eqnarray*}
M(t) = M(0) + \int_0^t\varphi(s)'dW(s), \enskip 0 \leq t \leq T, \textrm{ and } \int_0^T|\varphi(t)|^2dt <\infty \enskip \textrm{a.s.}
\end{eqnarray*}
\emph{In particular, $M$ has continuous sample paths a.s.}

Considerable effort has in the literature been made in order to find explicit formulas for the integrand $\varphi$, i.e. in order to find constructive representations of martingales, mainly using Malliavin calculus, see e.g. \cite{karatzas1991extension,malliavin2015stochastic,nualart2006malliavin,rogers2000diffusions} and the references therein. The recently developed \emph{functional It\^o calculus} includes a new type of constructive representation of square integrable martingales due to Cont and Fourni{\'e} see e.g. \cite{ContBook,contpathspace,Cont2013,cont2016weak}. The main result of the present paper is an extension of this result to local martingales. 

The organization of the paper is as follows. Section 2 is based on  \cite{ContBook} and contains a brief and heuristic account of the relevant parts of functional It\^o calculus including the constructive martingale representation theorem for square integrable martingales. Section 3 contains the local martingale extension of this theorem and a simple example.

\section{Constructive representation of square integrable martingales}
Denote an $n$-dimensional sample path by $\omega$. Denote a sample path stopped at $t$ by $\omega_t$, i.e. let $\omega_t(s) = \omega(t \wedge s), 0\leq s \leq T$. Consider a real-valued functional of sample paths $F(t,\omega)$ which is \emph{non-anticipative} (essentially meaning that $F(t,\omega) =F(t,\omega_t)$). The \emph{horizontal derivative} at $(t,\omega)$ is defined by
\begin{eqnarray*}
\mathcal{D} F(t,\omega) = \lim_{h \searrow  0}\frac{F(t+h,\omega_t)-F(t,\omega_t)}{h}.
\end{eqnarray*}
The \emph{vertical derivative} at $(t,\omega)$ is defined by $\nabla_\omega F(t,\omega)  = (\partial_iF(t,\omega), i=1,...,n)'$, where 
\begin{eqnarray*}
\partial_iF(t,\omega) = \lim_{h\rightarrow 0} \frac{F(t,\omega_t+  he_i I_{[t,T]})-F(t,\omega_t)}{h}.
\end{eqnarray*}
Higher order vertical derivatives are obtained by vertically differentiating vertical derivatives.

One of the main results of functional It\^o calculus is the functional It\^o formula, which is just the standard It\^o formula with the usual time and space derivatives replaced by the horizontal and vertical derivatives. If the functional $F$ is sufficiently regular (regarding e.g. continuity and boundedness of its derivatives), which we write as $F\in \bbbc^{1,2}_b$, then the functional It\^o formula holds, see \cite[ch.~5,6]{ContBook}. We remark that \cite{levental2013simple} contains another version of this result.

Using the functional It\^o formula it easy to see that if $Z$ is a martingale satisfying
\begin{eqnarray} \label{c12funcLINDENSJO}
Z(t) = F(t,W_t) \enskip  dt\times d\bbbp\textrm{-a.e., with } F\in\bbbc_b^{1,2}, 
\end{eqnarray}
then, for every $t\in [0,T]$, 
\begin{eqnarray*}
Z(t) = Z(0) + \int_0^t\nabla_\omega F(s,W_s)'dW(s) \enskip  a.s.
\end{eqnarray*}
We may therefore define the vertical derivative with respect to the process $W$ of a martingale $Z$ satisfying (\ref{c12funcLINDENSJO}) as the $dt\times d\bbbp$-a.e. unique process $\nabla_WZ$ given by
\begin{eqnarray} \label{zeqLINDENSJO}
\nabla_WZ(t) = \nabla_\omega F(t,W_t), \enskip  0 \leq t \leq T.
\end{eqnarray} 
Let $\mathcal{C}_b^{1,2}(W)$ be the space of processes $Z$ which allow the representation in (\ref{c12funcLINDENSJO}). Let $\mathcal{L}^2(W)$ be the space of progressively measurable processes $\varphi$ satisfying the condition $E[\int_0^T\varphi(s)'\varphi(s)ds]<\infty$. Let $\mathcal{M}^2(W)$ be the space of square integrable martingales with initial value $0$. Let $D(W) = \mathcal{C}_b^{1,2}(W)  \cap  \mathcal{M}^2(W)$.

It can be shown that $\{\nabla_WZ: Z \in D(W)\}$ is dense in $\mathcal{L}^2(W)$ and that $D(W)$ is dense in $\mathcal{M}^2(W)$ \cite[ch.~7]{ContBook}. Using this it is possible to show that the vertical derivative operator $\nabla_W(\cdot)$ admits a unique extension to $\mathcal{M}^2(W)$, in the following sense: \emph{For $Y \in \mathcal{M}^2(W)$ the (weak) vertical derivative $\nabla_W Y$ is the unique element in $\mathcal{L}^2(W)$ satisfying
\begin{eqnarray} \label{verticaldersqLINDENSJO}
E[Y(T)Z(T)] = E\left[\int_0^T\nabla_WY(t)'\nabla_WZ(t)dt\right]
\end{eqnarray} 
for every $Z\in D(W)$, where $\nabla_WZ$ is defined in (\ref{zeqLINDENSJO})}. The constructive martingale representation theorem (\cite[ch.~7]{ContBook}) follows: 
\begin{theorem}  [Cont and Fourni{\'e}] \label{cont-fornLINDENSJO} For any square integrable martingale $Y$ relative to $(\bbbp,\underline{\mathcal F})$ and every $t\in [0,T]$,
\begin{eqnarray*}
Y(t)  = Y(0) + \int_0^t\nabla_{W}Y(s)'dW(s)\enskip a.s.
\end{eqnarray*}
\end{theorem}

\section{Constructive representation of local martingales}
This section contains an extension of the vertical derivative $\nabla_W(\cdot)$ and the constructive martingale representation in Theorem \ref{cont-fornLINDENSJO} to local martingales. Let $\mathcal{M}^{\textrm{\tiny loc}}(W)$ denote the space of local martingales relative to $(\bbbp,\underline{\mathcal F})$ with initial value zero and RCLL sample paths. In Theorem \ref{extderivativeLINDENSJO} we extend the vertical derivative to $\mathcal{M}^{\textrm{\tiny loc}}(W)$. Using this extension we can formulate the constructive martingale representation theorem also for local martingales, see Theorem \ref{main-resLINDENSJO}.

Before extending the definition of the vertical derivative to $\mathcal{M}^{\textrm{\tiny loc}}(W)$ we recall the definition of a local martingale.

\begin{definition} \label{localMGdefLINDENSJO} $M$ is said to be a local martingale if there exists a sequence of non-decreasing stopping times $\{\theta_n\}$ with $\lim_{n \rightarrow \infty } \theta_n= \infty$ a.s. such that the stopped local martingale $M(\cdot\wedge\theta_n)$ is a martingale for each $n\geq 1$.
\end{definition}
 
\begin{theorem} [Definition of $\nabla_W(\cdot)$ on $\mathcal{M}^{\textrm{\tiny loc}}(W)$]\label{extderivativeLINDENSJO}\quad
\begin{itemize} 
\item{There exists a progressively measurable $dt\times d\bbbp$-a.e. unique extension of the vertical derivative $\nabla_W(\cdot)$ from $\mathcal{M}^2(W)$ to $\mathcal{M}^{\textrm{\tiny loc}}(W)$, such that, for $M \in \mathcal{M}^{\textrm{\tiny loc}}(W)$,
\begin{equation}
\begin{array}{rcl} \label{def-resLINDENSJO} 
M(t) & =&  \int_0^t\nabla_WM(s)'dW(s), 0 \leq t \leq T, \enskip \mbox{and} \nonumber\\
&& \int_0^T|\nabla_WM(t)|^2dt <\infty \enskip \mbox{a.s.}
\end{array}
\end{equation}}
 
\item{Specifically, for $M \in \mathcal{M}^{\textrm{\tiny loc}}(W)$ the vertical derivative $\nabla_WM$ is defined as the progressively measurable  $dt\times d\bbbp$-a.e. unique process satisfying
\begin{eqnarray}\label{sadadsLINDENSJO}
\nabla_WM(t) = \lim_{n\rightarrow\infty}\nabla_WM_n(t) \enskip dt\times d\bbbp\textrm{-a.e.}
\end{eqnarray} 
where $\nabla_WM_n$ is the vertical derivative of $M_n:= M(\cdot \wedge \tau_n) \in \mathcal{M}^2(W)$ and $\tau_n$ is given by 
\begin{eqnarray} \label{tau-stopLINDENSJO}
\tau_n = \theta_n \wedge \inf\{s\in [0,T]:|M(s)| \geq n\}\wedge T
\end{eqnarray}
where $\{\theta_n\}$ is an arbitrary sequence of stopping times of the kind described in Definition \ref{localMGdefLINDENSJO}.}
\end{itemize}
\end{theorem}
\begin{remark} \label{unique-remarkLINDENSJO} Note that if $M$ in Theorem \ref{extderivativeLINDENSJO} satisfies 
\begin{eqnarray*} 
M(t) = \int_0^t\gamma(s)'dW(s), 0 \leq t\leq T \enskip \mbox {a.s.}
\end{eqnarray*}
for some process $\gamma$, then $\gamma = \nabla_WM$ $dt\times d\bbbp$-a.e. It follows that the extended vertical derivative $\nabla_WM$ defined in Theorem \ref{extderivativeLINDENSJO} does not depend (modulo possibly on a null set $dt\times d\bbbp$) on the particulars of the chosen stopping times $\{\theta_n\}$.
\end{remark}
\begin{proof} The martingale representation theorem implies that, for $M \in \mathcal{M}^{\textrm{\tiny loc}}(W)$, there exists a progressively measurable process $\varphi$ satisfying
\begin{eqnarray} \label{pf1LINDENSJO}
\enskip M(t) =\int_0^t\varphi(s)'dW(s), \enskip 0 \leq t \leq T, \textrm{ and } \int_0^T|\varphi(t)|^2dt <\infty \enskip \textrm{a.s.}
\end{eqnarray}
Therefore, if we can prove that
\begin{eqnarray} \label{3434LINDENSJO}
\lim_{n\rightarrow \infty} \nabla_WM_n(t) = \varphi(t)\enskip dt\times d\bbbp\textrm{-a.e.},
\end{eqnarray}
then it follows that there exists a progressively measurable process, denote it by $\nabla_WM$, which is $dt\times d\bbbp$-a.e. uniquely defined by (\ref{sadadsLINDENSJO}) and satisfies 
\begin{eqnarray*}
\nabla_WM(t) = \varphi(t)\enskip dt\times d\bbbp\mbox{-a.e.,}
\end{eqnarray*}
which in turn implies that the integrals of $\nabla_WM$ and $\varphi$ coincide in the way that (\ref{pf1LINDENSJO}) implies (\ref{def-resLINDENSJO}). All we have to do is therefore to prove that (\ref{3434LINDENSJO}) holds.

Let us recall some results about stopping times and martingales. The stopped local martingale $M(\cdot\wedge\theta_n)$ is a martingale for each $n$, by Definition \ref{localMGdefLINDENSJO}.  Stopped RCLL martingales are martingales. The minimum of two stopping times is a stopping time and the hitting time 
\begin{eqnarray*}
\inf\{s\in [0,T]:|M(s)| \geq n\}
\end{eqnarray*}
is, for each $n$, in the present setting, a stopping time. Using these results we obtain that
\begin{eqnarray*}
M(\cdot \wedge \theta_n \wedge \inf\{s\in [0,T]:|M(s)| \geq n\}\wedge T) = M(\cdot\wedge\tau_n)
\end{eqnarray*}
is a martingale, for each $n$. Moreover, $M$ is by the standard martingale representation result a.s. continuous. Hence, we may define a sequence of, a.s. continuous, martingales $\{M_n\}$ by 
\begin{eqnarray} \label{weweLINDENSJO}
M_n = M(\cdot \wedge \tau_n) = \int_0^{\cdot \wedge \tau_n}\varphi(s)'dW(s) \enskip  \textrm{a.s.}
\end{eqnarray}
where the last equality follows from (\ref{pf1LINDENSJO}). Now, use the definition of $\tau_n$ in (\ref{tau-stopLINDENSJO}) to see that
\begin{eqnarray*}
|M_n(t)| = \left|\int_0^{t \wedge \tau_n}\varphi(s)'dW(s)\right| \leq n \enskip  \textrm{a.s.}
\end{eqnarray*} 
for any $t$ and $n$, and that in particular $M_n$ is, for each $n$, a square integrable martingale. Moreover, (\ref{weweLINDENSJO}) implies that $M_n$ satisfies 
\begin{eqnarray} \label{rtertxxLINDENSJO}
M_n(t) = \int_0^tI_{\{s \leq \tau_n\}}\varphi(s)'dW(s), \enskip 0 \leq t \leq T \enskip  \textrm{a.s.}
\end{eqnarray}
Since each $M_n$ is a square integrable martingale we may use Theorem \ref{cont-fornLINDENSJO} on $M_n$, which together with (\ref{rtertxxLINDENSJO}) implies that 
\begin{eqnarray} \label{rtertLINDENSJO}
M_n(t) = \int_0^t\nabla_WM_n(s)'dW(s)  =  \int_0^tI_{\{s\leq \tau_n\}}\varphi(s)'dW(s), \enskip 0 \leq t \leq T \enskip  \textrm{a.s.}
\end{eqnarray}
where $\nabla_WM_n$ is the vertical derivative of $M_n$ with respect to $W$ (defined in (\ref{verticaldersqLINDENSJO})) and where we also used the continuity of the It\^o integrals.  The equality of the two It\^o integrals in (\ref{rtertLINDENSJO}) implies that 
\begin{eqnarray} \label{rter2LINDENSJO}   
\nabla_WM_n(t)= I_{\{t \leq \tau_n\}}\varphi(t) \enskip dt\times d\bbbp\textrm{-a.e.}
\end{eqnarray} 
The local martingale property of $M$ implies that $\lim_{n \rightarrow \infty } \theta_n= \infty\enskip a.s.$ Using this and the definition of $\tau_n$ in (\ref{tau-stopLINDENSJO}) we conclude that for almost every $\omega \in \Omega$ and each $t \in [0,T]$ there exists an $N(\omega,t)$ such that 
\begin{eqnarray} \label{rter3LINDENSJO}
n\geq N(\omega,t) \Rightarrow \sup_{0\leq s \leq t}|M(\omega,s)|\leq n \textrm{ and } t\leq \theta_n(\omega) \Rightarrow t \leq \tau_n(\omega).
\end{eqnarray}  
It follows from (\ref{rter2LINDENSJO}) and (\ref{rter3LINDENSJO}) that there exists an $N(\omega,t)$ such that 
\begin{eqnarray*}  
n\geq N(\omega,t) \Rightarrow \nabla_WM_n(\omega,t) = \varphi(\omega,t) \enskip dt\times d\bbbp\textrm{-a.e.}
\end{eqnarray*}
which means that (\ref{3434LINDENSJO}) holds.
\qed\end{proof}

If $M$ is a RCLL local martingale then $M-M(0) \in \mathcal{M}^{\textrm{\tiny loc}}(W)$, which implies that $\nabla_W(M-M(0))$ is defined in Theorem \ref{extderivativeLINDENSJO}. This observation allows us to extend the definition of the vertical derivative to RCLL local martingales not necessarily starting at zero in the following obvious way.

\begin{definition} \label{generarldefLINDENSJO} The vertical derivative of a local martingale $M$ relative to $(\bbbp,\underline{\mathcal F})$ with RCLL sample paths is defined as the progressively measurable $dt\times d\bbbp\textrm{-a.e.}$ unique process $\nabla_WM$ satisfying
\begin{eqnarray}  \label{asdwwLINDENSJO}
\nabla_WM(t) = \nabla_W(M-M(0))(t), \enskip 0 \leq t \leq T,
\end{eqnarray}  
where $\nabla_W(M-M(0))(t)$ is defined in Theorem \ref{extderivativeLINDENSJO}.
\end{definition} 
The following result is an immediate consequence of Theorem \ref{extderivativeLINDENSJO} and Definition \ref{generarldefLINDENSJO}.
\begin{theorem}  \label{main-resLINDENSJO} If $M$ is a local martingale relative to $(\bbbp,\underline{\mathcal F})$ with RCLL sample paths, then
\begin{eqnarray*}
M(t)&=&M(0) + \int_0^t\nabla_WM(s)'dW(s), 0 \leq t \leq T, \enskip \mbox{and}\\
&& \int_0^T|\nabla_WM(t)|^2dt<\infty \enskip a.s.,
\end{eqnarray*} 
where $\nabla_WM(s)$ is defined in Definition \ref{generarldefLINDENSJO}.
\end{theorem} 
 
Let us try to clarify the theory by studying a simple example. It is straightforward to extend the results above to the case when the Wiener process $W$ is replaced by an adapted process $X$ given by 
\begin{eqnarray}  \label{example1LINDENSJO}
X(t)=X(0) + \int_0^t\sigma(s)dW(s),
\end{eqnarray}
where $\sigma$ is a matrix-valued adapted process satisfying suitable assumptions, mainly invertibility, see also \cite{ContBook,Cont2013}. Thus, a local martingale $M$ can be represented as
\begin{eqnarray*}
M(t)-M(0)= \int_0^t\nabla_W M(s)'dW(s) = \int_0^t\nabla_X M(s)'dX(s),
\end{eqnarray*}
and the relationship between the vertical derivatives with respect to $W$ and $X$ is $\nabla_W M(t)' = (\nabla_X M(t)')\sigma(t)$, cf. (\ref{example1LINDENSJO}). As an example consider the one-dimensional case and let $X$ with $X(0)=0$ be given by (\ref{example1LINDENSJO}) under the assumption that $\sigma(s)$ is a deterministic function of time and let $M$ be given by $M(t)= F(t,X_t)$ where $F$ is the non-anticipative functional $F(t,\omega) = \omega^3(t)- 3\int_0^t\omega(s)\sigma^2(s)ds$, i.e. let $M$ be the local martingale defined by
\begin{eqnarray*}
M(t)= X^3(t)- 3\int_0^tX(s)\sigma^2(s)ds.
\end{eqnarray*}
In this case the vertical derivative simplifies to the standard derivative i.e. $\nabla F_\omega(t,\omega) = 3\omega^2(t)$, see also \cite{ContBook,Cont2013} (we remark that the horizontal derivative is $\mathcal{D} F(t,\omega) = -3\omega(t)\sigma^2(t)$). In this case, $\nabla_X M(t) = 3X^2(t)$ and  
\begin{eqnarray*}
M(t) = \int_0^t3X^2(s)dX(s) = 
\int_0^t3 
X^2(s) 
\sigma(s)dW(s), \nonumber
\end{eqnarray*}
which we remark is easily found using the standard It\^o formula. Note that this also means that 
$\nabla_W M(t) 
=  3X^2(t)\sigma(t)
= \nabla_X M(t)\sigma(t)$.

\paragraph{Concluding remarks} 
Many of the applications that rely on martingale representation are within mathematical finance. A particular application that may benefit from the local martingale extension of the present paper is optimal investment theory, in which the discounted (using the state price density) optimal wealth process is a (not necessarily square integrable) martingale, see e.g. \cite[ch.~3]{karatzas1}, see also \cite{lindensjo2016optimal}. In particular, using functional It\^o calculus it is possible to derive an explicit formula for the optimal portfolio in terms of the vertical derivative of the discounted optimal wealth process, see also \cite{lindensjo-portfolio}. Similar explicit formulas for optimal portfolios based on the Malliavin calculus approach to constructive martingale representation have, under restrictive assumptions, been studied extensively, see e.g. \cite{benth2003explicit,detemple2005closed,di2009optimal,lakner1998optimal,lakner2006portfolio,ocone1991generalized,okur2010white,pham2001optimal}. The general connection between Malliavin calculus and functional It\^o calculus is studied in e.g. \cite{ContBook,Cont2013}.

\begin{acknowledgement}
The author is grateful to Mathias Lindholm for helpful discussions.
\end{acknowledgement}


%
%
%

\end{document}